\documentclass [12 pt, letterpaper] {article}

        \usepackage {times,amssymb,amsmath,array,calc,enumerate}
   {\begin {list}{} %{\textbf{Example}\arabic{example}.}
   {\usecounter{example}%
   \small%
   \setlength{\labelsep}{0pt}\setlength{\leftmargin}{40pt}%
   \setlength{\labelwidth}{0pt}%
   \setlength{\listparindent}{0pt}}}%
   {\end{list}}

\renewcommand {\theequation} {\@arabic\c@equation}

\makeatother
   {\begin {list}{} %{\textbf{Example}\arabic{example}.}
   {%\usecounter{example}
   \small%
    \setlength{\itemindent}{10pt}
    \setlength{\labelsep}{5pt}\setlength{\leftmargin}{40pt}%
    \setlength{\labelwidth}{5pt}%
   \setlength{\itemsep}{-1.8pt}
    \setlength{\listparindent}{0pt}
   }
   }%
   {\end{list}}

\newcommand{\N}{\mathcal {N}}

\newcommand{\A}{\mathcal{A}}
\newcommand{\M}{\mathcal{M}}

\renewcommand {\theequation} {\arabic{equation}}
\newtheorem{Def}{Definition}[section]
\newtheorem{Th}{Theorem}[section]
\newtheorem{Col}{Corollary}[section]
\newtheorem{Prop}{Proposition}[section]
\newtheorem{lemma}{Lemma}[section]

\begin {document}

\centerline{\Large {\bf  Type II$_1$ Factors With A Single
Generator}}

\thispagestyle{empty}   $\bigskip$

\centerline{\large  Junhao \  Shen}

\bigskip

\centerline{Mathematics Department, University of New Hampshire,
Durham, NH, 03824}

\bigskip

\centerline{email:  \qquad jog2@cisunix.unh.edu}

$\bigskip$

\noindent\textbf{Abstract: } In the paper, we study the generator
problem of type II$_1$ factors. By defining a new concept closely
related to the numbers of generators of a von Neumann algebra, we
are able to show that a large class of type II$_1$ factors are
singly generated, i.e., generated by two self-adjoint elements. In
particular, we show that most of type II$_1$ factors,
  whose free entropy dimensions   are known to be less than or equal to one, are singly generated.

\vspace{0.2cm} \noindent{\bf Keywords:} generator problem, type
II$_1$ factor, free entropy dimension

\vspace{0.2cm} \noindent{\bf 2000 Mathematics Subject
Classification:} Primary 46L10, Secondary 46L54

\section{Introduction}

 Let $H$ be a separable complex Hilbert space, $\mathcal B(H)$ be
  the algebra consisting of all bounded linear operators from
  $H$ to $H$. A von Neumann algebra $\mathcal M$ is defined to be a
   self-adjoint subalgebra of $\mathcal B(H)$ which is closed in the
   strong
  operator topology. Factors are the von Neumann algebras whose
  centers are scalar multiples of the identity.  The factors
   are  classified by means of a relative dimension
  function into type I, II, III factors. (see [9])

  The generator problem for von Neumann algebras is the question of
   whether every von Neumann
  algebra
  acting on a separable Hilbert space
   can be
  generated by two self-adjoint elements    (equivalently be  singly generated). It is a
  long-standing open problem (see [8]), and is still unsolved. Many
  people (see [3], [4], [6], [9], [12], [13], [14], [19] )
  have contributed  to this topic. For example, von Neumann   [10] proved  that every abelian von Neumann algebra
   is generated by one self-adjoint element and every type II$_1$
   hyperfinite von
   Neumann algebra
   is singly generated. W. Wogen    [19] showed that every properly infinite von Neumann
   algebra is singly generated.  It follows   that
     the generator problem for   von Neumann algebras, except for the
   non hyperfinite type II$_1$ von Neumann algebras, is solved.
   In [14],  S. Popa proved that a type II$_1$ von Neumann algebra
    with a Cartan
   subalgebra is singly generated.
  L. Ge and S. Popa   [6] proved that certain type II$_1$ factors
   are   singly generated. These type II$_1$ factors  include
  the ones
   with property $\Gamma$, those that are not  prime. In [4], Ge and the
   author proved that some type II$_1$ factors with property T, including
   $   L(SL(\Bbb Z, 2m+1))$ ($m\ge 1$), are singly generated. This
   result answered one question proposed by Voiculescu.

    In the early 1980s, D.
Voiculescu began the development of  the theory of free probability
and free entropy. This new and powerful tool was crucial in solving
some  old open problems in the field of von Neumann algebras. In his
ground-breaking paper [16], Voiculescu defined a new concept called
``free entropy dimension", by which Voiculescu was able to show that
free group factors have no Cartan subalgebras [17]. To better
understand the free entropy dimension of von Neumann algebras has
become an urgent task for the subject. It is believed that the free
entropy dimension is closely related to the number of generators of
a von Neumann algebra. Also because Voiculescu's free entropy
dimension is always defined on finitely generated von Neumann
algebras, it is important to know whether a   von Neumann algebra is
finitely generated,  even singly generated, or not.

Inspired by Voiculescu's approach to free entropy dimension in [17],
in this paper  we define a new concept, $\mathcal G(\M)$, of a
diffuse von Neumann algebra $\M$. More specifically, suppose $\M$ is
a diffuse von Neumann algebra with a tracial state $\tau$. The
invariant $\mathcal G(\M)$ is designed to count the numbers of
generators when these generators are in matricial forms (see
definitions in section 2). $\mathcal G(\M)$ has many good
properties, some of which are listed as follow.
\begin{enumerate}
 \item If $\M$ is a type II$_1$ factor and $\mathcal G(\M) < 1/4$, then
 $\M$ is singly generated.
 \item If $\M$ is a diffuse hyperfinite von Neumann algebra with a tracial state $\tau$, then
 $\mathcal G(\M)=0$.
 \item Suppose that $\M$ is a type II$_1$ factor with the tracial state $\tau$.
Suppose $\{\N_k\}_{k=1}^\infty$ is a sequence of von Neumann
subalgebras of $\M$   that generates $\M$ as a von Neumann algebra
and   $\N_k\cap \N_{k+1}$ is a diffuse von Neumann subalgebra of
$\M$ for each $k\ge 1$.  If  $\mathcal G(\N_k)=0$ for $k\ge 1$, then
$ \mathcal G(\M)= 0.$ In particular, $\M$ is singly generated.
\item Suppose that $\M$ is a type II$_1$ factor with the tracial state $\tau$.
 Suppose $\N$ is a von Neumann subalgebra of $\M$ with $\mathcal G(\N)=0$ and $\{u_k\}_{k=1}^\infty$ is a
 family of unitary elements in $\M$ such that $\{\N, u_1,u_2,\ldots
 \ \}$ generates $\M$ as a von Neumann algebra. If there exists a
 family of Haar unitary elements $\{v_{k,n}\}_{k,n=1}^\infty $ in $\N$ such that
$\lim_{n\rightarrow \infty} dist_{\|\cdot
\|_2}(u_k^*v_{k,n}u_k,\N)=0$ for $k\ge 1  ,$    then $\mathcal
G(\M)=0$. In particular, if   there exists a
 family of Haar unitary elements $\{v_{k } \}_{k= 1}^\infty $ in $\N$ such that
$  u_k^*v_{k }u_k $ is contained in $\N$ for $k\ge 1  ,$    then
$\mathcal G(\M)=0$. And $\M$ is singly generated.
\item  Suppose that $\M$ is a type II$_1$ factor with the tracial state $\tau$.
 Suppose $\N$ is a von Neumann
subalgebra of $\M$ with $ \mathcal G(\N)=0$ and
$\{u_k\}_{k=1}^\infty$ is a
 family of Haar unitary elements in $\M$ such that $\{\N, u_1,u_2,\ldots
 \ \}$ generates $\M$ as a von Neumann algebra. If $u_1$ is contained in
 $\N$ and $u_{k+1}^*u_ku_{k+1}^*$ is contained in the von Neumann subalgebra
 generated by $\N\cup\{ u_1, \ldots, u_k\}$ for $k\ge 1$,  then $\mathcal
 G(\M)=0$. In particular, $\M$ is singly generated.
\end{enumerate}
Using the listed properties of $\mathcal G(\M)$, we can easily
provide many new examples of singly generated type II$_1$ factors
besides the ones we have known; and give  the new proofs of some
known results, such as   type II$_1$ factors with property $\Gamma$
or the ones with Cartan subalgebras, nonprime type II$_1$ factors,
some type II$_1$ factors with property T, are singly generated type
II$_1$ factors. More importantly, the class of type II$_1$ factors
considered in [5] and [7] is a class of singly generated type II$_1$
factors. Thus most of the type II$_1$ factors, whose free entropy
dimensions are known to be lass than or equal to one, are singly
generated.

\vspace{0.2cm}

Let us fix our notation. For a subset  $A$ of $ \mathcal B( H)$, let
$A''$ denote the von Neumann algebra generated by the elements of
$A$ in $\mathcal B(H)$.

\section{Definitions}

Suppose that $\M$ is a diffuse von Neumann algebra  with a tracial
state $\tau$.

\begin{Def}
Let $\{p_j\}_{j=1}^k$ be a family of mutually orthogonal projections
of $\M$ with $\tau(p_j)=1/k$ for each $1\le j\le k$. For an element
$x$ of $\M$, we define
$$
\mathcal I(x;\{p_{ j}\}_{ j=1}^k) =\frac {\left |\{ (i,j) \ | \
p_{i}xp_{j}\ne 0  \}\right |}{k^2},
$$ where $|   \cdot   |$ denotes the cardinality of the set; and the
support of $x$ on $\{p_j\}_{j=1}^k$ is defined by
$$
\mathcal S(x;\{p_{ j}\}_{ j=1}^k) = \vee \{p_j\ | \ p_jx\ne 0, \text
{ or } \ xp_j \ne 0, \ 1\le j\le k\},
$$ where $\vee$ denotes the union of the projections.
 For elements
$x_1,\ldots, x_n$ in $\M$, we define
$$
\mathcal I(x_1,\ldots,x_n;\{p_j\}_{j=1}^k) =\sum_{m=1}^n \mathcal
I(x_m;\{p_j\}_{j=1}^k).$$
\end{Def}

\begin{Def}

For each positive integer $k$, let $\frak E_k$ denote the collection
of  all $\{p_{j}\}_{j=1}^k$, the    families of mutually orthogonal
   projections of $\M$ with $\tau(p_j)=1/k$ for each $1\le
j\le k$. Suppose $x_1,\ldots,x_n$ are elements in $\M$. We define
$$
\mathcal I (x_1,\ldots,x_n; k) = \inf \ \{ \mathcal
I(x_1,\ldots,x_n;\{p_{j}\}_{j=1}^k )\ | \ \ \{p_{j}\}_{j=1}^k \in
\frak E_k\};
$$and
$$
\mathcal G(\M; k) = \left \{\begin{aligned} &\inf \ \{\mathcal I
(x_1,\ldots,x_n; k) \ | \ \ \ \text {$x_1,\ldots,x_n$
generate $\M$ as a von Neumann algebra.}\ \}\\ &\quad \\
 & \infty; \ \ \text {if $\M$ is not finitely generated.} \end{aligned}\right .
$$Then, we define
$$
\mathcal G(\M)\ \ = \ \ \liminf_{k\rightarrow \infty} \ \mathcal
G(\M; k).
$$
\end{Def}

\noindent {\bf Remark: } By the definition, for every $k>1$, we know
that $ \mathcal G(\M; k^n)$ is a decreasing function as $n$
increases. Thus, $\mathcal G(\M)\le  \mathcal G(\M; k)\le \mathcal
G(x_1,\ldots, x_n;k)$ for each $k\ge 1$, each family of generators
$\{x_1,\ldots,x_n\}$ of $\M$.

\section{$\mathcal G(\M): \M $  is  a Diffuse Hyperfinite von Neumann Algebra   }
In this section, we are going to compute $\mathcal G(\M)$ when $\M$
is a diffuse hyperfinite von Neumann algebra with a tracial state
$\tau$.

\begin{lemma} Suppose  $\M=\M_1\oplus \M_2$ is a von Neumann algebra with
a tracial state $\tau$, where $\M_1,\M_2$ are the von Neumann
subalgebras of $\M$. Then $\mathcal G(\M)\le  \mathcal
G(\M_1)+\mathcal G(\M_2) $.
 \end{lemma}
 \noindent{\bf Proof: }
 It is trivial when one of $\mathcal
G(\M_1),\mathcal G(\M_2)$ is infinite. Let $c_i=\mathcal G(\M_i)$
for $i=1,2$.
  By the definitions of $ \mathcal
G(\M_1)$ and $\mathcal G(\M_2) $, for each positive $\epsilon,$ we
know there
 exist a large positive integer $k$,
elements $\{p_{j }\}_{j =1}^{  k}$, $\{q_{j }\}_{j =1}^{  k}$,
$\{x_1,\ldots, x_n\}$ and $\{y_1,\ldots,y_m\}$ of $\M$ such that
\begin{enumerate}
\item $\{p_{j }\}_{j =1}^{  k}$, or
$\{q_{j }\}_{j =1}^{  k}$, is a family of mutually orthogonal
   projections of $\M_1$, or $\M_2$ respectively, with
$\tau(p_{j })=\tau(I_{\M_1})/k, \tau(q_{j })=\tau(I_{\M_2})/k$,
$\sum_{j }p_{j }=I_{\M_1}$ and $\sum_{j }q_{j }=I_{\M_2}$.
\item $\{x_1,\ldots, x_n\}$, or
$\{y_1,\ldots,y_m\}$, is a family of generators of $\M_1$, or $\M_2$
respectively.
\item $$\begin{aligned}
  &\mathcal
   I(x_1,\ldots,x_n;\{p_{j }\}_{j =1}^{  k})\le
  c_1+\epsilon\\
  &\mathcal I(y_1,\ldots,y_m;\{q_{j }\}_{j =1}^{  k })\le
  c_2+\epsilon.
\end{aligned}$$
\end{enumerate}
A little computation shows
$$
\mathcal I(x_1,\ldots,x_n,y_1,\ldots,y_m;\{p_{j}+ q_{j}\}_{j=1 }^{
k})\ \le { c_1+c_2+2\epsilon  } .
$$Hence, by definitions, we have $\mathcal G(\M)\le
  c_1+c_2 +2\epsilon$; whence $\mathcal G(\M)\le \mathcal
G(\M_1)+\mathcal G(\M_2) $.

\vspace{0.2cm}

The following two propositions are obvious.
\begin{Prop}
 Suppose $M_k$ is a factor of type I$_k$ and $\{e_{ij}\}_{i,j=1}^k$
 is a system of matrix units of $M_k$. Let $x_1= e_{11}$ and $x_2=
 \sum_{i=1}^{k-1} (e_{i,i+1}+e_{i,i+1}^*).$ Then $x_1, x_2$ are two
 self-adjoint elements that generate $M_k$ as a von Neumann
 algebra.
\end{Prop}
\begin{Prop}
Suppose $\M\simeq \mathcal A\otimes \N$ is a von Neumann algebra
with a tracial state $\tau$, where $\A,\N$ are finitely generated
von Neumann subalgebras of $\M$. If $\A$ is a von Neumann subalgebra
with $\mathcal G(\A)=0$, then $\mathcal G (\M)=0$. In particular,
if $\A$ is a diffuse abelian von Neumann subalgebra of $\M$, then
$\mathcal G(\M)=0$.
\end{Prop}
\begin{Th}
Suppose that $\mathcal R$ is the hyperfinite type II$_1$ factor.
Then $\mathcal G (\mathcal R) =0$.
\end{Th}
\noindent {\bf Proof: } Let $\{n_k\}_{k=1}^\infty$ be a sequence of
positive integers with $n_k\ge 3$ for $k=1,2,\ldots $ . It is easy
to see that $\mathcal R \simeq \otimes_{k=1}^\infty  M_{n_k}(\Bbb
C)$ where $  M_{n_k}(\Bbb C)$ is the algebra of $n_k\times n_k$
matrices with complex entries. Assume that
$\{e_{i,j}^{(k)}\}_{i,j=1}^{n_k}$ is the canonical system of matrix
units of $  M_{n_k}(\Bbb C)$. We should identify  $  M_{n_k}(\Bbb
C)$ with its canonical image in $\otimes_{k=1}^\infty  M_{n_k}(\Bbb
C)$ if it causes no confusion. Let
$$
  \begin{aligned}
     x_1 & =  e_{11}^{(1)} + \sum_{k=1}^\infty \frac 1 {2^k}
     e_{22}^{(1)}\otimes e_{22}^{(2)}\otimes \cdots \otimes
     e_{22}^{(k)}\otimes e_{11}^{(k+1)}\\
     x_2  &= \sum_{j=2}^{n_1} (e_{j-1,j}+e_{j,j-1}) +
     \sum_{k=1}^\infty\sum_{j=2}^{n_{k+1} } \frac 1 {2^k} e_{22}^{(1)}\otimes e_{22}^{(2)}\otimes \cdots \otimes
     e_{22}^{(k)}  \otimes (e_{j-1,j}^{(k+1)}+e_{j,j-1}^{(k+1)})
   \end{aligned}
$$
It is easy to check that, for each $k\ge 1$,
$\{e_{ij}^{(k)}\}_{i,j=1}^{n_k}$ is in the von Neumann subalgebra
generated by $x_1,x_2$ in $\mathcal R$. Thus $x_1, x_2$ are two
self-adjoint elements  that generate $\mathcal R$. Moreover, we have
$\mathcal I(x_1,x_2; \{e_{jj}^{(1)}\}_{ j=1}^{n_1})\le 3/n_1.$
Therefore, $\mathcal G(\mathcal R)\le 3/n_1$, for all $n_1\ge 3$.
Hence $\mathcal G(\mathcal R)=0.$ Q.E.D.

\vspace{0.2cm}

 Now we are able to compute
$\mathcal G(\M)$ for a  diffuse hyperfinite von Neumann algebra
$\M$.
 \begin{Th}
Suppose $\M$ is a diffuse hyperfinite von Neumann algebra with a
tracial state $\tau$. Then $\mathcal G(\M)=0$.

 \end{Th}
 \noindent {\bf Proof: } Note a diffuse hyperfinite von Neumann algebra $\M$
 can always be decomposed as
 $$
\M\simeq \A_0\otimes \mathcal R \oplus \left ( \oplus_{k=1}^\infty
\A_k\otimes \M_{n_k}(\Bbb C) \right ),
 $$ where $\mathcal R$ is the hyperfinite type II$_1$ factor, $\A_0$ is an abelian von Neumann
subalgebra of $\M$,
  and $\A_k$
 is a diffuse abelian von Neuamnn subalgebra of $\M$. The rest follows from Lemma 3.1, Proposition 3.1 and 3.2, and
 Theorem 3.1. \qquad Q.E.D.

 \vspace{0.2cm}

\section{Cut-and-Paste theorem}

The proof of following easy theorem is based on a ``cut-and-paste"
trick from [6] or [4].
 \begin{Th}
Suppose that $\M$ is a von Neumann algebra with a tracial state
$\tau$. Suppose $\{e_{ij}\}_{i,j=1}^k$ is a system of matrix units
of a I$_k$ subfactor of $\M$ with $\sum_{j=1}^k e_{jj}=I$. If
$x_1,\ldots, x_n$ are the elements in $\M$ such that $$\mathcal
I(x_1,\ldots,x_n;\{e_{jj}\}_{j=1}^k)=c^2
$$ with $c\le \frac 1 2 - \frac 1 { k},$ then there exists a projection
$q$ in $\{x_1,\ldots, x_n,e_{ij}; 1\le i,j\le k\}''$ so that
$$\tau(\mathcal S(q;\{e_{jj}\}_{j=1}^k )) \le 2c+ \frac 2 { k}$$ and
\qquad \qquad\quad\quad $\{q, e_{ij}; 1\le i,j\le k\}''$
$=\{x_1,\ldots, x_n,e_{ij}; 1\le i,j\le k\}''.$
\end{Th}
\noindent {\bf Proof: } Let $$ \mathcal T = \{(i,j,p)\ |  \
e_{ii}x_pe_{jj}\ne 0, \  \ 1\le i,j\le k, 1\le p\le n\}.
$$ Note that
$$\begin{aligned} |\mathcal T|=k^2\cdot \mathcal
 I(x_1,\ldots,x_n;\{e_{jj}\}_{j=1}^k)=c^2k^2,
\end{aligned}
$$ and the cardinality of the set
$$
\{  (s,t) \ |  \  1\le s\le [ck]+1, \  [ck]+2\le t\le 2[ck]+2\}.
$$ is equal to $([ck]+1)^2\ge c^2k^2$.
There exists an injective mapping from  $( i,j,p)\in \mathcal T $ to
$$
(s,t)\in \{  (s,t) \ |  \  1\le s\le [ck]+1, \  [ck]+2\le t\le
2[ck]+2\},
$$ and denote this map by $( i,j,p) \mapsto
(s(i,j,p),t(i,j,p))$. Then each $e_{ii}x_pe_{jj}$ may be replaced by
$e_{s(i,j,p)i}x_pe_{jt(i,j,p)}$ for all $( i,j,p) \in \mathcal T$.
Let
$$\begin{aligned}
 y &= \sum_{( i,j,p) \in\mathcal T} \ \ \left (e_{s(i,j,p)i}x_pe_{jt(i,j,p)}+ \left
(e_{s(i,j,p)i}x_pe_{jt(i,j,p)}\right )^*\right )\\
 q_1&=\sum_{s=1}^{[ck]+1} e_{ss} \qquad \text { and } \qquad q_2= \sum_{t=[ck]+2}^{2[ck]+2} e_{tt}
.\end{aligned}
$$ Without loss of generality, we can assume that $\|y\|\le 1$. Then
let
$$
 q= \frac 1 2 q_1(1+(1-y^2)^{1/2})q_1 +\frac 1 2
y+ \frac 1 2 q_2(1-(1-y^2)^{1/2})q_2.
$$
By the construction of $q$, it is easy to check that $q$ is a
projection in $\M$ with $\tau(\mathcal S(q;\{e_{jj}\}_{j=1}^k )) \le
2c+ \frac 2 { k}$ and  $\{q, e_{ij}; 1\le i,j\le k\}''
=\{x_1,\ldots, x_n,e_{ij}; 1\le i,j\le k\}''.$ \qquad  Q.E.D.

\vspace{0.2cm} The following theorem indicates the relationship
between $\mathcal G(\M)$ and singly generated type II$_1$ factors.
\begin{Th}
Suppose $\M$ is a type II$_1$ factor with the tracial state $\tau$.
If $\mathcal G(\M) < 1/4$, then $\M$ is singly  generated.
\end{Th}
\noindent {\bf Proof: } Note that $\M$ is a type II$_1$ factor. From
the preceding theorem and the definition of $\mathcal G(\M)$, for a
sufficiently large integer $k$, there exist a system of matrix
units, $\{e_{ij}\}_{i,j=1}^k$, of a I$_k$ subfactor of $\M$ and a
projection $q$ in $\M$ so that the following hold. (i) $\{q\}\cup
\{e_{ij}\}_{ i,j=1}^k$ generates $\M$; and (ii) $\tau(\mathcal S(q ;
\{e_{jj}\}_{ j=1}^k) )\le 2\sqrt{\mathcal G(\M)}+ \frac 2 { k} < 1 -
\frac 1 k.$ Therefore we can   assume that $e_{11}$ and $q$ are two
mutually orthogonal projections of $\M$. Let $x_1= e_{11} + 2q$ and
$x_2=
 \sum_{i=1}^{k-1} (e_{i,i+1}+e_{i,i+1}^*).$ By same argument as in Proposition 3.1 and
 the fact that $\{q\}\cup
\{e_{ij}\}_{ i,j=1}^k$ generates $\M$, we obtain that
 $x_1, x_2$ are two self-adjoint elements of $\M$ that generate $\M$
 as a von Neumann algebra.\qquad Q.E.D.

 \vspace{0.2cm}

 \noindent{\bf Remark: } Instead of constructing a projection $q$ in
 Theorem 4.1, if  we are interested in constructing a self-adjoint element, then the result in
 Theorem 4.2 can be improved as follows. {\em
Suppose $\M$ is a type II$_1$ factor with the tracial state $\tau$.
If $\mathcal G(\M) < 1 /\sqrt 2$, then $\M$ is singly  generated.}

\section{Main Results}

\vspace{0.2cm}The following lemma essentially comes from Popa's
remarkable paper [15].
\begin{lemma}
Suppose $\M$ is a type II$_1$ factor with the tracial state $\tau$.
Suppose $\{p_j\}_{j=1}^k$ is a family of mutually orthogonal
   projections in $\M$ with each $\tau(p_j)=1/k$. Then there
exists a hyperfinite type II$_1$ subfactor $\mathcal R$  of $\M$
such that $\mathcal R'\cap M=\Bbb C I$ and $\{p_{j}\}_{j=1}^k
\subset \mathcal R$.

\end{lemma}
\noindent{\bf Proof: } By [15], there exists a  hyperfinite
subfactor $\mathcal R_0$ of $\M$ such that $\mathcal R_0'\cap
\M=\Bbb C I$. Since $\M$ is a type II$_1$ factor, there exists a
unitary element $w$ in $\M$ such that $\{p_{j}\}_{j=1}^k \subset
w^*\mathcal R_0w$. Let $\mathcal R= w^*\mathcal R_0w$. Then
$\mathcal R$ is a hyperfinite type II$_1$ subfactor of $\M$ such
that $\mathcal R'\cap \M=\Bbb C I$ and $\{p_{j}\}_{j=1}^k \subset
\mathcal R$.

\vspace{0.2cm}
 \noindent {\bf Definition: } \ $\M_1$ \ is called an
irreducible subfactor of a type II$_1$ factor \ $\M$ \ if \
$\M_1\subset\M$ \ and \ $\M_1'\cap \M=\Bbb C I$.

\begin{lemma}
Suppose that $\M$ is a type II$_1$ factor with the tracial state
$\tau$. Suppose $\N$ is a von Neumann subalgebra of $\M$ with
$\mathcal G(\N)=c$. Then for each $\epsilon>0$, there exists an
irreducible subfactor $\M_\epsilon$ of $\M$ such that \ $\N\subset
\M_\epsilon\subset \M$ and \ $\mathcal G(\M_\epsilon)\le
c+\epsilon.$

\end{lemma}
\noindent {\bf Proof: } Since $\mathcal G(\N)=c$, there exist some
positive integer $k> \frac 8 \epsilon$, a family of mutually
orthogonal    projections $\{p_j\}_{j=1}^k$ in $\N$ with
$\tau(p_j)=1/k$ for $1\le j\le k$, and a family of generators
$\{x_1,\ldots, x_n\}$ of $\N$, such that
$$
\mathcal I(x_1, \ldots,x_n; \{p_j\}_{j=1}^k)\le c+\frac   \epsilon
2.
$$
By Lemma 5.1, we can find an irreducible hyperfinite type II$_1$
subfactor $\mathcal R$ of $\M$ such that $\{p_j\}_{j=1}^k\subset
\mathcal R$. Thus there exists a system of matrix units
 $\{e_{ij}\}_{i,j=1}^k$  of a I$_k$ subfactor  $M_k$  of $\mathcal
R$ such that $e_{jj}=p_j$ for each $j=1,\ldots,k$. It is easy to see
that $\mathcal R\simeq \mathcal R_1\otimes M_k$ for some hyperfinite
type II$_1$ subfactor $\mathcal R_1$ of $\mathcal R$. By Theorem 3.1
and 4.2, we know the hyperfinite subfactor $\mathcal R_1$ is
generated by two self-adjoint elements $y_1, y_2$ that commute with
$M_k$. By Proposition 3.1, $M_k$ is generated by two self-adjoint
elements $z_1=p_1$ and $z_2=\sum_{j=1}^{k-1} (e_{j,j+1}+ e_{j+1,j})$
as a von Neumann algebra. A little computation shows that
$$
\mathcal I(x_1, \ldots,x_n, y_1,y_2, z_1,z_2; \{p_j\}_{j=1}^k)\le
c+\frac \epsilon 2 + \frac 2 k + \frac 2 k\le c+ \epsilon.
$$ Let $\M_\epsilon$ be the von Neumann subalgebra generated by $\mathcal R$
and $\N$ in $\M$. Since $\mathcal R$ is an irreducible type II$_1$
subfactor of $\M$, $\M_\epsilon$ is also an irreducible type II$_1$
subfactor of $\M$, which is   generated by $x_1,\ldots,
x_n,y_1,y_2,z_1,z_2$ in $\M$ as a von Neumann algebra. Hence
$\mathcal G(\M_\epsilon)\le c+\epsilon.$ \quad  Q.E.D.

\vspace{0.2cm}

\noindent {\bf Definition: } The family of elements
$\{e_{ij}\}_{i,j=1}^k$ is called a subsystem of matrix units of a
von Neumann algebra $\M$ if the following hold: (i)
$\{e_{ij}\}_{i,j=1}^k\subset \M$; (ii) there exists a projection $p$
in $\M$ such that $\sum_{j=1}^k e_{jj}=p$; (iii) $e_{ij}^*=e_{ji}$
for $1\le i,j\le k$; (iv) $e_{il}e_{lj}=e_{ij}$ for $1\le i,l,j\le
k.$

\vspace{0.2cm}

Next proposition is our main technical result in the paper.

\begin{Prop}
Suppose that $\M$ is a type II$_1$ factor with the tracial state
$\tau$. Suppose $\{\N_k\}_{k=1}^\infty$ is a sequence of von Neumann
subalgebras of $\M$   such that $\{\N_k\}_{k=1}^\infty$ generates
$\M$ as a von Neumann algebra and  $\N_k\cap \N_{k+1}$ is a diffuse
von Neumann subalgebra of $\M$ for all $k\ge 1$. Suppose, for each
$k\ge 1$, $\epsilon>0$, there is an irreducible subfactor \
$\M_{k,\epsilon}$ \ of $\M$ such that \ $\N_k\subset
\M_{k,\epsilon}\subset \M$ \ and \ $\mathcal G(\M_{k,\epsilon})\le
\epsilon.$ Then \ $ \mathcal G(\M)= 0.$ In particular, $\M$ is
singly generated.
\end{Prop}
\noindent{\bf Proof: } Let $\epsilon<1/8 $ be a positive number.
From assumption on $\N_1$, there exists an irreducible type II$_1$
subfactor $\M_1$ of $\M$ such that $\N_1\subset \M_1\subset \M$ and
$\mathcal G(\M_1)\le \epsilon.$
   By Theorem 4.1 and the
definition of $\mathcal G(\M_1)$, for a sufficiently large integer
$m_1>\frac 3 { \epsilon}$, there exist a projection $q_1$ in
$\mathcal M_1$ and a system of matrix units
$\{e_{ij}^{(1)}\}_{i,j=1}^{m_1}$ of  $\mathcal M_1$ such that
$\sum_{j=1}^k e_{jj}^{(1)}=I$, $\tau (\mathcal
S(q_1;\{e_{jj}^{(1)}\}_{j=1}^{m_1}) )\le 3\epsilon$, and
$\{q_1\}\cup \{e_{ij}^{(1)}\}_{i,j=1}^{m_1}$ generates $\mathcal
M_1$ as a von Neumann algebra. Without loss of generality, we can
assume that $e_{11}^{(1)}, e_{22}^{(1)}, q$ are mutually orthogonal
projections in $\mathcal M_1$.

We claim that we can construct a sequence of positive integers
$\{m_k\}_{k=1}^\infty$,  a sequence of irreducible type II$_1$
subfactors $\M_k$ of $\M$, subsystems of matrix units
$\{\{e_{ij}^{(k)}\}_{i,j=1}^{m_k}\}_{k=1}^\infty$, and  a family of
projections $\{q_k\}_{k=1}^\infty$, such that  \   (i)   \
$\N_k\subset \M_k\subset \M$ for $k\ge 1$; \ (ii) \
$\sum_{j=1}^{m_{k+1}} e_{jj}^{(k+1)}= e_{22}^{(k)}$ for $k\ge 1$; \
(iii) $q_{k+1}=e_{22}^{(k)}q_{k+1}e_{22}^{(k)}$, \
$q_{k+1}e_{11}^{(k+1)}=0$, \ $q_{k+1}e_{22}^{(k+1)}=0$  for $k\ge
1$; (iv) $\{\M_1,\ldots\ ,\M_k\}''$ $=\{ q_1, \ldots, q_k,
e_{ij}^{(p)}; \
 1\le i,j\le  m_p, 1\le p\le k \}''$ for $k\ge 1$.

We have already finished the construction when $k=1$. Suppose that
 we have finished the   construction  till $k$-step.
 Note that,  by the assumption on $\N_{k+1}$, there  exists an irreducible subfactor $\M_{k+1}$ of $\M$ such that
$$\mathcal G(\M_{k+1})\le \left (\frac 1{ 8m_1\cdots m_k }\right
)^2, $$and $\N_{k+1}\subset \M_{k+1}\subset \M$, i.e., {\em (i)
holds.}

 By the definition of $\mathcal G(\M_{k+1})$,
 there exist a sufficiently large   integer $m_{k+1}>10
 $, a family of mutually  orthogonal    projections
 $\{p_{j}\}_{j=1}^{m_1\cdots m_{k+1}}$
 in $\M_{k+1}$    with  each   $\tau(p_j)=1/m_1\cdots m_{k+1}$ and a family of
 generators
 $\{x_1,\ldots, x_n\}$ of  $\mathcal M_{k+1}$ such that
 \begin{align}
\mathcal I(x_1,\ldots,x_n;\{p_{j}\}_{j=1}^{m_1\cdots
m_{k+1}})\le\left (\frac 1{ 4m_1\cdots m_k }\right )^2\tag{$\ast$}
 \end{align}
From the induction hypothesis on each $\M_j$, we know that
$\{\M_j\}_{j=1}^k$ are a family of irreducible type II$_1$
subfactors of $\M$, which implies $\{\M_1,\ldots,\M_k\}''$ is a type
II$_1$ subfactor of $\M$. And,
\begin{align}\{ q_1, \ldots, q_k,  e_{ij}^{(p)}; \
 1\le i,j\le  m_p, 1\le p\le k \}''=\{\M_1,\ldots,\M_k\}''.\tag{$\ast$$\ast$}\end{align} Let
$\{e_{ij}^{({k+1})}\}_{i,j=1}^{m_{k+1}}$ be a subsystem of matrix
units in   $\{\M_1,\ldots,\M_k\}''$ such that
$e_{22}^{(k)}=\sum_{j=1}^{m_{k+1}}e_{jj}^{({k+1})}$, i.e., {\em (ii)
holds.}

Then $$ \mathcal T_{k+1}=\{e_{i_12}^{(1)}\cdots
e_{i_{k},2}^{(k)}e_{st}^{({k+1})} e_{2,j_{k}}^{(k)}\cdots
e_{2j_1}^{(1)} \ | \ 1\le i_p,j_p\le m_p, 1\le p\le k, 1\le s,t\le
m_{k+1} \}$$ is a system of matrix units of a I$_{m_1m_2\cdots
m_km_{k+1}}$ subfactor of $\{\M_1,\ldots,\M_k \}''$; and
$$\mathcal P_{k+1}=\{e_{i_12}^{(1)}\cdots e_{i_{k},2}^{(k)}e_{ss}^{({k+1})}
e_{2,i_{k}}^{(k)}\cdots e_{2i_1}^{(1)} \ | \ 1\le i_p \le m_p, 1\le
p\le k, 1\le s \le m_{k+1}\}$$ is a family of mutually orthogonal
equivalent projections in $\{\M_1,\ldots,\M_k \}''$ with sum $I_{\M}
$. Note the following facts: (1) $\M_k\cap \M_{k+1}$ is a diffuse
von Neumann subalgebra; (2) $\mathcal P_{k+1}$ is in the type II$_1$
subfactor  $\{\M_1,\ldots, \M_k\}''$; (3)
$\{p_{j}\}_{j=1}^{m_1\cdots m_{k+1}}$ is in the type II$_1$
subfactor $\M_{k+1}$. Thus there exist unitary elements $v_{k+1}$ in
$\{\M_1,\ldots,\M_k\}''$ and $w_{k+1}$ in $\M_{k+1}$ such that
$w_{k+1}v_{k+1}$  maps $\mathcal P_{k+1}$, one to one, onto
$\{p_j\}_{j=1}^{m_1\cdots m_{k+1}}$. By ($\ast$), it is easy to
compute that
$$
\mathcal I( v_{k+1}^*w_{k+1}^*x_1w_{k+1}v_{k+1}, \ldots,
 v_{k+1}^*w_{k+1}^*x_nw_{k+1}v_{k+1};\mathcal P_{k+1} ) \le\left (\frac 1{ 4m_1\cdots m_k }\right )^2
$$ By Theorem 4.1, there exists a projection $q_{k+1}$ in $\M$ so that
\begin{align}\{v_{k+1}^*w_{k+1}^*x_1w_{k+1}v_{k+1},\quad  \ldots,\quad
 v_{k+1}^*w_{k+1}^*x_nw_{k+1}v_{k+1}, \mathcal T_{k+1} \}''= \{q_{k+1}, \mathcal T_{k+1}\}''\tag{$\ast$$\ast$$\ast$}\end{align}  and $$\tau(\mathcal
S(q_{k+1};\mathcal P_{k+1} )) \le \frac 1{ 2m_1\cdots m_k }+ \frac
2{ m_1\cdots m_{k+1}}<\frac 1 {m_1\cdots m_k}- \frac 3 {m_1\ldots
m_{k+1}}.$$ Because
$$
\tau(e_{22}^{(k)})=\frac 1 {m_1\cdots m_k}\ , \qquad
\tau(e_{11}^{(k+1)})=\tau(e_{22}^{(k+1)})=\frac 1 {m_1\cdots
m_km_{k+1}} \ ,
$$
 we might assume that $q_{k+1}=e_{22}^{(k)}q_{k+1}e_{22}^{(k)}$, \
$q_{k+1}e_{11}^{({k+1})}=0$, \ $q_{k+1}e_{22}^{({k+1})}=0$, i.e.,
 {\em (iii) holds.}

Note $v_{k+1}$ is  in $\{\M_1,\ldots,\M_k\}''$, which, by
($\ast$$\ast$), is in the von Neumann algebra generated by $\{ q_1,
\ldots, q_k\} \cup\{ \{e_{ij}^{(p)}\}_{i,j=1,\ldots, m_p; 1\le p\le
k}\}$. It is easy to check that
$$
\{\mathcal T_{k+1}\}'' = \{ e_{ij}^{(p)}; \
 1\le i,j\le  m_p, 1\le p\le k+1 \}''.
$$
 Together with ($\ast$$\ast$$\ast$), we get that $\{ w_{k+1}^*x_1w_{k+1} ,\quad \ldots,\quad
 w_{k+1}^*x_nw_{k+1}  \}$ is contained in the von Neumann subalgebra generated by $\{
q_1, \ldots, q_k,q_{k+1}\} \cup\mathcal T_{k+1}$ in $\M$. But
$\M_{k+1}$ is generated by $\{ w_{k+1}^*x_1w_{k+1} ,\quad
\ldots,\quad
 w_{k+1}^*x_nw_{k+1}  \}$, since $w_{k+1}$ is in $\M_{k+1}$. Hence, $\M_{k+1}$ is
   in the von Neumann algebra generated by $\{
q_1, \ldots, q_k,q_{k+1}\} \cup\mathcal T_{k+1}$. Combining with the
facts that $$\begin{aligned} \{\M_1,\ldots\ ,\M_k\}''  =\{  q_1,
\ldots, q_k, e_{ij}^{(p)}; \
 1\le i,j\le  m_p, 1\le p\le k \}''  &\supset \mathcal T_{k+1},\\
 q_{k+1}   \in \{v_{k+1}^*w_{k+1}^*x_1w_{k+1}v_{k+1},\
\ldots,\
 v_{k+1}^*w_{k+1}^*x_nw_{k+1}v_{k+1}, \mathcal T_{k+1} \}''&\subset \{\M_1,\ldots, \M_{k+1},\mathcal T_{k+1}  \}'',\end{aligned}$$ we
know that
$$\begin{aligned}
 \{\M_1,&\ldots\ ,\M_k,\M_{k+1}\}''  \subset
\{ \M_{k+1}, q_1, \ldots, q_k,   e_{ij}^{(p)}; \
 1\le i,j\le  m_p, 1\le p\le k+1 \}''\\
&\subset \{ \M_{k+1}, q_1, \ldots,  q_{k },\mathcal T_{k+1}\}''
 \subset \{  q_1, \ldots,  q_{k+1 },\mathcal T_{k+1}\}''\\
 &\subset \{  q_1, \ldots  ,q_{k+1}, e_{ij}^{(p)}; \
 1\le i,j\le  m_p, 1\le p\le k +1 \}''\\
&\subset\{\M_1,\ldots\ ,\M_k,\M_{k+1},\mathcal
T_{k+1}\}''\subset\{\M_1,\ldots\ ,\M_k,\M_{k+1}\}'';
\end{aligned}  $$
   whence {\em (iv) holds.}  This finishes the construction at $(k+1)$-th step.

Let
$$
  \begin{aligned}
     x_1 &=  \left (\sum_{k=1}^\infty \frac 1 {2^k} e_{11}^{(k)}\right
     ) + \left (  \sum_{k=1}^\infty   \frac 1 {3^k}  q_k  \right )\\
     x_2 &  =
     \sum_{k=1}^\infty\sum_{j=2}^{m_k } \frac 1 {2^k}
    (e_{j-1,j}^{(k)}+e_{j,j-1}^{(k)})
  \end{aligned}
$$
Note that, by induction hypothesis (iii), we know $\{e_{11}^{(k)},
q_k; k\ge 1\}$ is a family of mutually orthogonal projections in
$\M$. Thus, $\{e_{11}^{(k)}, q_k; k\ge 1\}$  is in the von Neumann
subalgebra generated by $x_1$. By the construction of $x_2$ and the
fact that $\{e_{11}^{(k)};k\ge 1\}$ is in the von Neumann subalgebra
generated by $x_1$, we get that $\{e_{ij}^{(k)}; \
 1\le i,j\le  m_k, k\ge 1 \}$ is in the von Neumann subalgebra
 generated by $\{x_1, x_2\}$. Hence, by induction hypothesis (iv), $\{\M_k\}_{k=1}^\infty$ is in the von Neumann subalgebra
 generated by $\{x_1, x_2\}$, i.e.,  $x_1, x_2$ are self-adjoint elements in
$\mathcal M$ that generate $\mathcal M$ as a von Neumann algebra.
Moreover, a little computation shows that
$$
\mathcal I(x_1, x_2; \{e_{ij}^{(1)}\}_{i,j=1}^{m_1}) \le 3\epsilon+
\frac 3 {m_1}\le 4\epsilon.
$$ Therefore, $\mathcal G(\mathcal M)\le 4\epsilon,$ for all
$\epsilon>0$. It follows that $\mathcal G(\mathcal M)=0$.\quad
Q.E.D.

\vspace{0.2cm}

\begin{Th} Suppose that $\M$ is a type II$_1$ factor with the tracial state $\tau$.
Suppose $\{\N_k\}_{k=1}^\infty$ is a sequence of von Neumann
subalgebras of $\M$   that generates $\M$ as a von Neumann algebra
and   $\N_k\cap \N_{k+1}$ is a diffuse von Neumann subalgebra of
$\M$ for each $k\ge 1$.  If  $\mathcal G(\N_k)=0$ for $k\ge 1$, then
$ \mathcal G(\M)= 0.$ In particular, $\M$ is singly generated.
\end{Th}
\noindent {\bf Proof: }  The result follows easily from Lemma 5.2
and Proposition 5.1.

\vspace{0.2cm}

\noindent{\bf Definition: } Suppose that $\M$ is a diffuse von
Neumann subalgebra with a tracial state $\tau$. A unitary element
$v$ in $\M$ is call a Haar unitary element if $\tau(v^m)=0$ for all
$m\ne 0$.

\begin{lemma}
Suppose that $\M$ is a type II$_1$ factor with the tracial state
$\tau$. Suppose $\N$ is a von Neumann subalgebra of $\M$ such that
$\mathcal G(\N)=c$. Suppose $u$ is a unitary element in $\M$ such
that, for some Haar unitary element $v$ in $\N$, $u^*vu$ is
contained in $\N$. Then, for every $\epsilon>0$, there exists an
irreducible type II$_1$ subfactor $\M_\epsilon$ such that \
$\{\N\cup\{u\}\}''\subset \M_\epsilon\subset \M$ \ and \ $\mathcal
G(\M_\epsilon)\le c+\epsilon.$
\end{lemma}
\noindent {\bf Proof: } By Lemma 5.2,   there exists an irreducible
type II$_1$ subfactor $\N_\epsilon$ of $\M$ such that $\N\subset
\N_\epsilon\subset \M$ and $\mathcal G(\N_\epsilon)\le
c+\epsilon/2.$ Thus, by the definition of $\mathcal G(\N_\epsilon)$,
there exist some positive integer $k> \frac 8 \epsilon$, a family of
mutually orthogonal    projections $\{p_j\}_{j=1}^k$ in
$\N_\epsilon$ with $\tau(p_j)=1/k$  for $1\le j\le k$, and a family
of generators $\{x_1,\ldots, x_n\}$ of $\N_\epsilon$, such that
$$
\mathcal I(x_1, \ldots,x_n; \{p_j\}_{j=1}^k)\le c+ \frac  \epsilon
2.
$$
Note $u$ is a unitary element in $\M$ such that, for some Haar
unitary element $v$ in $\N$, $u^*vu$ is contained in $\N$. It
follows that there exist two families of mutually orthogonal
   projections, $\{e_{j}\}_{j=1}^k $, $\{f_{j}\}_{j=1}^k$,
in $\N$ with $\tau(e_j)=\tau(f_j)=1/k$ such that $u^*e_ju=f_j$ for
$j=1, \ldots,k$. Note $\N_\epsilon$ is a type II$_1$ subfactor that
contains $\N$. There exist two unitary elements $w_1, w_2$ in
$\N_\epsilon$ such that $p_j=w_1^*e_jw_1=w_2^*f_jw_2$ for
$j=1,\ldots,k$. Thus $w_1^*uw_2p_j=p_jw_1^*uw_2$ for $j=1,\ldots,k$.
It follows
$$
\mathcal I(x_1, \ldots,x_n,w_1^*uw_2; \{p_j\}_{j=1}^k)\le c+\frac
\epsilon 2+ \frac 1 k \le c+\epsilon.
$$Let $\M_\epsilon$ be the von Neumann subalgebra generated by $x_1, \ldots,x_n,w_1^*uw_2$
in $\M$; whence $\mathcal G(\M_\epsilon)\le c+ \epsilon.$ Note
$\N_\epsilon$ is contained in $\M_\epsilon$, so are $w_1,w_2$. Thus
$u$ is also contained in $\M_\epsilon$, whence
$\{\N\cup\{u\}\}''\subset \M_\epsilon\subset \M$.  From the fact
that $\N_\epsilon'\cap \M=\Bbb CI$, it follows that $\M_\epsilon$ is
an irreducible type II$_1$ subfactor of $\M$. Q.E.D.

\vspace{0.2cm}

\begin{Th}
 Suppose that $\M$ is a type II$_1$ factor with the tracial state $\tau$.
 Suppose $\N$ is a von Neumann subalgebra of $\M$ and $\{u_k\}$ is a
 family of unitary elements in $\M$ such that $\{\N, u_1,u_2,\ldots
 \ \}$ generates $\M$ as a von Neumann algebra. Suppose there exists a
 family of Haar unitary elements $\{v_k\}_{k=1}^\infty $  such that $u_k^*v_ku_k$
 is in $\N$ for $k\ge 1$. If $ \mathcal G(\N)=0, $ then $\mathcal
 G(\M)=0$. In particular, $\M$ is singly generated.
\end{Th}
\noindent{\bf Proof: } Let   $\N_k$ be the von Neumann subalgebra
generated by $\N$ and $u_k$ in $\M$ for $k\ge 1$. Using Lemma 5.3
and Proposition 5.1, we  easily obtained the result. Q.E.D.
\begin{Th}
 Suppose that $\M$ is a type II$_1$ factor with the tracial state $\tau$.
 Suppose $\N$ is a von Neumann subalgebra of $\M$ and $\{u_k\}$ is a
 family of Haar unitary elements in $\M$ such that $\{\N, u_1,u_2,\ldots
 \ \}$ generates $\M$ as a von Neumann algebra. Suppose $u_1$ is in
 $\N$ and $u_{k+1}^*u_ku_{k+1}^*$ is in the von Neumann subalgebra
 generated by $\N\cup\{ u_1, \ldots, u_k\}$ for $k\ge 1$. If $ \mathcal G(\N)=0, $ then $\mathcal
 G(\M)=0$. In particular, $\M$ is singly generated.
\end{Th}
\noindent{\bf Proof: } Let  $\N_k$ be the von Neumann subalgebra
generated by $\N$ and $u_1, \ldots,u_k $ in $\M$ for $k\ge 1$. Using
Lemma 5.3, inductively, and Proposition 5.1, we can easily obtained
the result. Q.E.D.
\begin{lemma}
Suppose that $\M$ is a type II$_1$ factor with the tracial state
$\tau$. Suppose $\N$ is a von Neumann subalgebra of $\M$ such that
$\mathcal G(\N)=0$. Suppose $u$ is a unitary element in $\M$ such
that, for a family of   Haar unitary elements $\{v_n\}$ in $\N$,
$\lim_{n\rightarrow \infty} dist_{\|\cdot \|_2}(u^*v_nu,\N)=0$.
Then, for every $\epsilon>0$, there exists an irreducible type
II$_1$ factor $\M_\epsilon$ such that $\{\N\cup\{u\}\}''\subset
\M_\epsilon\subset \M$ and $\mathcal G(\M_\epsilon)\le \epsilon.$
\end{lemma}
\noindent{\bf Proof: } By the assumption on the unitary element $u$,
we can assume that there exists a family of Haar unitary elements
$\{w_n\}_{n=1}^\infty$ in $\N$ such that $\lim_{n\rightarrow
0}\|u^*v_nu-w_n\|_2=0$ (see [7]). Equivalently, for every $k>1$,
$\epsilon=1/k$, there exist $\{p_j\}_{j=1}^k$, $\{q_j\}_{j=1}^k$
families of mutually orthogonal    projections of $\N$ with each
$\tau(p_j)=\tau(q_j)=1/k$, such that $\|u-\sum_{j=1}^k
p_juq_j\|_2<\epsilon$. Let $x_k=\sum_{j=1}^k p_juq_j$ and
$\N_k=\{\N,x_k\}''$. Thus $x_k\stackrel{\|\cdot\|_2}\longrightarrow
u$. A straightforward adaption of the proofs of  Lemma 5.3 and
Proposition 5.1 shows that there exist a subsequence
$\{k_p\}_{p=1}^\infty$ of $\{k\}_{k=1}^\infty$ and an irreducible
subfactor $\M_\epsilon$ of $\M$ such that $\{\N_{k_p}\}_{p=1}^\infty
\subset \M_\epsilon\subset \M$ and $\mathcal G(\M_\epsilon)\le
\epsilon.$ But $x_{k_p}\in \N_{k_p}$ and $x_{k_p}
\stackrel{\|\cdot\|_2}\longrightarrow u$, as $p\rightarrow \infty$.
Thus $u\in \M_\epsilon$. This completes the proof. \quad Q.E.D.

\vspace{0.2cm} \noindent{\bf Remark: } The statement that ``$u$ is a
unitary element in $\M$ such that, for a family of   Haar unitary
elements $\{v_n\}$ in $\N$, $\lim_{n\rightarrow \infty}
dist_{\|\cdot \|_2}(u^*v_nu,\N)=0$" is equivalent to say that ``the
constant sequence $
 (u)_n $ is a unitary element in $\M^\omega$ such that, for some Haar
unitary elements $(v_n)_n,(w_n)_n$ in $ \N^\omega$,
$[(u)_n]^*[(v_n)_n][(u)_n]=[(w_n)_n]$," \  where $\omega$ is a free
ultra-filter of $\Bbb N$ and $\M^\omega$, or $\N^\omega$, is the
corresponding ultra-power of $\M$, or $\N$ respectively, along
$\omega$.

\vspace{0.2cm}
 Using Lemma 5.4 and Proposition 5.1, we can easily obtain the following
theorem.
\begin{Th}
 Suppose that $\M$ is a type II$_1$ factor with the tracial state $\tau$.
 Suppose $\N$ is a von Neumann subalgebra of $\M$ and $\{u_k\}$ is a
 family of unitary elements in $\M$ such that $\{\N, u_1,u_2,\ldots
 \ \}$ generates $\M$ as a von Neumann algebra. Suppose there exists a
 family of Haar unitary elements $\{v_{k,n}\}_{k,n=1}^\infty $  in $\N$  such that
$\lim_{n\rightarrow \infty} dist_{\|\cdot
\|_2}(u_k^*v_{k,n}u_k,\N)=0$ for $k\ge 1  .$ If $ \mathcal G(\N)=0,
$ then $\mathcal G(\M)=0$. In particular, $\M$ is singly generated.
\end{Th}

\section{A Few Applications}
In this section, we show that many type II$_1$ factors   are singly
generated. We provide new proofs of some known results, such as type
II$_1$ factors with property $\Gamma$ or the ones with Cartan
subalgebras are singly generated. New examples of singly generated
type II$_1$ factors are also given.

\vspace{0.2cm}

\noindent Using Theorem 5.3 in [1] (or Lemma 4 in [7]) and Theorem
3.2  and 5.4, we have the following result from [6].
\begin{Col}
Suppose $\M$ is a type II$_1$ factor with property $\Gamma$. Then
$\mathcal G(\M)=0$. In particular, $\M$ is singly generated.
\end{Col}
 \noindent  Using Theorem 3.2  and 5.3, we have
the following result.
\begin{Th}
Suppose that $\M$ is a type II$_1$ factor. Suppose that
$\{u_k\}_{k=1}^\infty$ is a family of Haar unitary elements in $\M$
that generate $\M$ and $u_{k+1}^*u_ku_{k+1}$ is contained in the von
Neumann subalgebra generated by $\{u_1,\ldots,u_k\}$ for $k\ge 1$.
Then $\mathcal G(\M)=0$. In particular, $\M$ is singly generated.
\end{Th}
As another corollary of  Theorem 3.2  and 5.3, we obtain the
following result from [14].
\begin{Col}
Suppose $\M$ is a type II$_1$ factor with Cartan subalgebras. Then
$\mathcal G(\M)=0$. In particular, $\M$ is singly generated.
\end{Col}
Using Theorem 3.2,  5.2 and 5.3, we have the following result.
\begin{Th}
 Suppose $\mathcal{M}$ is a   type II$_{1}$ factor
generated by a family $\left\{  u_{ij} \right\}_{i,j=1}^\infty  $ of
Haar unitary elements in $\mathcal{M}$ such that
\begin{enumerate}
\item for each $i,j$, $u_{i+1,j}^*u_{ij}u_{i+1,j} $ is in the von Neumann
subalgebra generated by $\{u_{1j},\ldots,u_{ij}\};$

\item for each $j\geq1$, $\left\{  u_{1j},u_{2j},\ldots\right\}
\bigcap\left\{  u_{1,j+1},u_{2,j+1},\ldots\right\} \neq\varnothing.$
\end{enumerate}
Then $\mathcal G(\M)=0$. In particular, $\M$ is singly generated.
\end{Th}
 The following corollaries follow easily from Theorem 6.2 (also see [4], [6]).
\begin{Col}
Suppose $\M=L(SL(\Bbb Z, 2m+1))$ $(m\ge 1)$ is the group von Neumann
algebra associated with $SL(\Bbb Z, 2m+1)$, the special linear group
with integer entries. Then $\mathcal G(\M)=0$. In particular, $\M$
is singly generated.
\end{Col}
\begin{Col}
Suppose $\mathcal M$ is a nonprime type II$_1$ factor, i.e.
$\mathcal M\simeq \mathcal N_1\otimes \mathcal N_2$ for some type
II$_1$ subfactors $ \mathcal N_1, \mathcal N_2$ of $\M$. Then
$\mathcal G(\M)=0$. In particular, $\M$ is singly generated.
\end{Col}

\noindent{\bf Remark: } Combining with the results in [5], [7], we
have shown that most of the type II$_1$ factors, whose free entropy
dimensions are known to be less than or equal to one, are singly
generated.

\vspace{0.2cm}

\noindent{\bf Examples:} New examples of singly generated II$_1$
factors can be constructed by considering the group von Neumann
algebras associated with some countable discrete groups. The
following are a few of them. (i) Let $G$ be the group $\langle g_1,
g_2,\ldots \ | \ g_ig_{i+1}=g_{i+1}g_i, i=1,\ldots \rangle$. Then
$\mathcal G(L(G))=0$ and $L(G)$ is singly generated, where $L(G)$ is
the group von Neumann algebra associated with $G$. \ (ii) Let $G$ be
the group $\langle a,b,c \ |  \ ab^2a^{-1}=b^3,
ac^2a^{-1}=c^3\rangle$. Then $\mathcal G(L(G))=0$ and $L(G)$ is
singly generated.

\vspace{1cm} \small{}

\end{document}